\documentclass{amsproc}

\usepackage{amssymb}

\usepackage{amsmath}

\newtheorem{thm}{Theorem}[section]

\newtheorem{cor}[thm]{Corollary}

\newtheorem{con}[thm]{Conjecture}

\theoremstyle{definition}

\newtheorem{example}[thm]{Example}

\theoremstyle{remark}

\numberwithin{equation}{section}

\newcommand\C{{\mathbb{C}}}

\newcommand\Q{{\mathbb Q}}

\newcommand\N{{\mathbb{N}}}

\newcommand\ps{{\rm ps}}

\newcommand\expon{{\rm ex}}

\newcommand\bq{\begin{equation}}
\newcommand\eq{\end{equation}}
\newcommand\beq{\begin{eqnarray*}}
\newcommand\eeq{\end{eqnarray*}}
\newcommand\ben{\begin{enumerate}}
\newcommand\een{\end{enumerate}}
\newcommand\bit{\begin{itemize}}
\newcommand\eit{\end{itemize}}

\newcommand\des{{\rm des}}
\newcommand\exc{{\rm exc}}
\newcommand\inv{{\rm inv}}
\newcommand\maj{{\rm maj}}

\newcommand\sg{{\mathfrak S}}
\newcommand\Des{{\rm DES}}
\newcommand\Dex{{\rm DEX}}

\newcommand\ch{{\rm ch}}

\newcommand\x{{\mathbf x}}

\newcommand\Exp{{\rm Exp}}

\newcommand\asc{{\rm asc}}

\newcommand\sym{\mathfrak S}

\def\h{{\mathsf h}}

\def\ff{{\mathbb F}}
\def\cc{{\mathbb C}}
\def\pp{{\mathbb P}}

\def\inc{{\rm{inc}}}

\def\hess{{\mathcal H}}

\def\Ps{{\rm \bf ps}}

\def\ee{{\mathsf e}}
\def\p{{\mathsf p}}

\def\q{\bold q}

\newlength\cellsize 
\savebox2{%
\begin{picture}(15,15)
\put(0,0){\line(1,0){15}}
\put(0,0){\line(0,1){15}}
\put(15,0){\line(0,1){15}}
\put(0,15){\line(1,0){15}}
\end{picture}}
\newcommand\cellify[1]{\def\thearg{#1}\def\nothing{}%
\ifx\thearg\nothing
\vrule width0pt height\cellsize depth0pt\else
\hbox to 0pt{\usebox2\hss}\fi%
\vbox to 15\unitlength{
\vss
\hbox to 15\unitlength{\hss$#1$\hss}
\vss}}
\newcommand\tableau[1]{\vtop{\let\\=\cr
\setlength\baselineskip{-16000pt}
\setlength\lineskiplimit{16000pt}
\setlength\lineskip{0pt}
\halign{&\cellify{##}\cr#1\crcr}}}
\savebox3{%
\begin{picture}(15,15)
\put(0,0){\line(1,0){15}}
\put(0,0){\line(0,1){15}}
\put(15,0){\line(0,1){15}}
\put(0,15){\line(1,0){15}}
\end{picture}}
\newcommand\expath[1]{%
\hbox to 0pt{\usebox3\hss}%
\vbox to 15\unitlength{
\vss
\hbox to 15\unitlength{\hss$#1$\hss}
\vss}}

\begin{document}

\title[From Poset Topology ... Chromatic Symmetric Functions]{From Poset Topology to $q$-Eulerian Polynomials to Stanley's Chromatic
Symmetric Functions}
\author[Shareshian]{John Shareshian$^1$}
\address{Department of Mathematics, Washington University, St. Louis, MO 63130}
\thanks{$^{1}$Supported in part by NSF Grant
 DMS 1202337}
\email{shareshi@math.wustl.edu}

\author[Wachs]{Michelle L. Wachs$^2$}
\address{Department of Mathematics, University of Miami, Coral Gables, FL 33124}
\email{wachs@math.miami.edu}
\thanks{$^{2}$Supported in part by NSF Grant
DMS 1202755}

\subjclass[2010]{}

\date{August 31, 2014}

\begin{abstract} In recent years we have worked on a project involving poset topology, various analogues of Eulerian polynomials, and a refinement of Richard Stanley's chromatic symmetric function.  Here we discuss how Stanley's ideas and results have influenced and inspired our own work.
\end{abstract}

\maketitle

\section{A walk in the woods}  It is a privilege and an honor to contribute an article to a volume celebrating Richard Stanley's 70th birthday.  Over the years, Richard's work has had a 
 tremendous influence on the research of so many of us working in algebraic and geometric combinatorics.  It is our pleasure to discuss in this article the impact of Richard's work on our  ongoing project involving Eulerian polynomials and symmetric functions.   Working on this project has been very much like taking a walk in a beautiful forest.  At every turn, we meet Richard Stanley, and each time we run into him our walk gets much more interesting.

The study of connections between  poset topology  and permutation enumeration has its roots in Stanley's pioneering work on edge labelings of posets,  Cohen-Macaulay complexes, and  group actions on posets.   In Section~\ref{posetsec}  we describe how Stanley's work on this topic inspired our work on the topology of a $q$-analog of a certain poset introduced by Bj\"orner and Welker in connection with commutative algebra.  We also describe how our work on this poset led to the discovery of  a new $q$-analog, involving the major index and excedance number, of  Euler's exponential generating function formula for the Eulerian polynomials.

We begin Section~\ref{q-EulerSec} by presenting a $q$-analog of Euler's exponential generating function formula due to Stanley.  Stanley's formula involves the inversion number and the descent number. 
 Next, we describe our own $q$-analog, which is obtained from a symmetric function identity using  stable principal specialization.  
 After discussing our symmetric function identity, we present in the remainder of Section~\ref{q-EulerSec} some consequences and refinements of our main results.   One of these consequences is the  $q$-unimodality of  the $q$-Eulerian polynomials that appear in our $q$-analog.   Here we find a close connection to earlier work of Stanley, in which the hard Lefschetz theorem is used to prove unimodality results in combinatorics.

 In Section~\ref{csfsec}, we we explain how the symmetric function arising in the previous sections is an instance of a refinement of Stanley's chromatic symmetric function for graphs.  We were inspired to define this refinement by an observation of Stanley about our earlier work on $q$-Eulerian polynomials discussed in Section \ref{q-EulerSec}.  The refinement is a polynomial in $t$, whose coefficients are in general  not symmetric functions, but rather  quasisymmetric functions.  However, if the graph at hand is the incomparability graph $G$  of a unit interval order  (with an appropriate labeling of the vertices), the coefficients are symmetric.  
 
 We present in the remainder of Section~\ref{csfsec} various results and conjectures on our chromatic quasisymmetric functions.  Of particular interest is a conjecture that asserts that with $G$ as above, the chromatic quasisymmetric function of $G$ is $\ee$-positive and $\ee$-unimodal.  This reduces to the Stanley-Stembridge $\ee$-positivity conjecture for unit interval orders when $t$ is set equal to 1.  We describe how our $\ee$-unimodality conjecture and an exercise in Stanley's EC1 on generalized Eulerian polynomials led  to a conjecture connecting the chromatic quasisymmetric functions to  a certain class of subvarieties of flag varieties due to De Mari and Shayman, known as the regular semisimple Hessenberg varieties.  
 Our conjecture asserts that, with $G$ as above,  the chromatic quasisymmetric function of  $G$ is essentially the Frobenius characteristic of the representation of the symmetric group on the cohomology of the Hessenberg variety naturally associated to $G$. This conjecture  implies Schur-positivity (which we prove by other means) of the chromatic quasisymmetric function of $G$.  Schur-unimodality (which is still open) will follow from our conjecture and the hard Lefschetz theorem.   We hope that our conjecture will lead eventually to a proof of the Stanley-Stembridge $\ee$-positivity conjecture.

  We assume that the reader is familiar with various ideas from algebraic, enumerative and topological combinatorics.  Luckily and unsurprisingly, all of these ideas are explained in the books \cite{St53,StEC1,StEC2}.

\section{Poset topology} \label{posetsec} Our project began while we were participating in a combinatorics program organized by Anders Bj\"orner and Richard Stanley, at  the Mittag-Leffler Institute.  
We worked on a conjecture of  Bj\"orner and Welker  \cite{BjWe}
on poset topology and permutation enumeration. 
Interest in this topic can be traced  back to Stanley's seminal 
work  on edge labelings of posets in \cite{St9, St11,St21} and on Cohen-Macaulay complexes  in \cite{St27,St34,St36,St52,St53}.   The fact that Cohen-Macaulay posets in general, and such posets admitting a group action in particular, can be useful in enumeration problems is made clear in \cite{St50}.  

We give a brief description of some work of Stanley on the rank-selected M\"obius invariant, which  inspired the poset topology part of our project.  For a ranked and bounded poset $P$ of length $n$ with rank function $r_P$, and a subset $S \subseteq [n-1]$, we consider the rank selected subposet $P_S:=\{x \in P : r_P(x) \in S\} \cup \{\hat 0, \hat 1\}.$  In \cite{St9} Stanley showed that if $P$ is a  distributive lattice then there is a way to label the edges of the Hasse diagram of $P$ so that the rank-selected M\"obius invariant 
$\mu_{P_S}(\hat 0, \hat 1) $ is equal to $(-1)^n$ times the number maximal chains whose label sequences have descent set $S$. (Throughout this paper, the M\"obius function on a poset $P$ will be denoted by $\mu_P$.)  This result is generalized to the class of supersolvable lattices in \cite{St11} and further  generalized to the class of admissible lattices in \cite{St21}.  
For the Boolean algebra $B_n$ on $n$ elements, the label sequences in \cite{St9} correspond to  the permutations in the symmetric group $\sg_n$.  Thus  for each $S \subseteq [n-1]$, the rank-selected M\"obius invariant is given by 
\begin{equation} \label{rankseleq} \mu_{(B_n)_S}(\hat 0, \hat 1) = (-1)^n | \{ \sigma \in \sg_n : \Des(\sigma) = S\}|,\end{equation}  where  $\Des(\sigma)$ is the descent set of $\sigma$.

 The  theory of lexicographic shellability, which relates edge labelings of posets to the topology of  order complexes of the posets,  
 has its roots in Stanley's theory of 
admissible labelings.  Indeed,    Stanley conjectured, and Bj\"orner proved in \cite{Bj} that all admissible lattices are Cohen-Macaulay. This means that the homology of the order complex of each open interval in such a lattice vanishes below the top dimension.  In proving this conjecture,   Bj\"orner was led  to introduce the theory of lexicographical shellability in \cite{Bj}.   This theory, which was further developed by Bj\"orner and Wachs in \cite{BjWa1,BjWa2,BjWa3,BjWa4}, has proved to be an important tool  for establishing Cohen-Macaulayness of posets and  determining the homotopy type of  order complexes; see also \cite{Wa}.

Through his work in  \cite{ St11, St21, St28, St50},  
Stanley demonstrates that  one can obtain beautiful  $q$-analogs of results  involving  the 
M\"obius function of the Boolean algebra   by replacing $B_n$ with the lattice $B_n(q)$ of subspaces of the $n$-dimensional vector space over the finite field $\ff_q$.  For example,
the following $q$-analog of (\ref{rankseleq}) was obtained by Stanley in \cite[Theorem 3.1]{St28}.  
 For each $S \subseteq [n-1]$,
\begin{equation} \label{qrankseleq}\mu_{(B_n(q))_S}(\hat 0, \hat 1) = (-1)^n \sum_{\substack {\sigma\in \sg_n \\ \Des(\sigma) = S}} q^{\inv(\sigma)},\end{equation}
where $\inv$ is the number of inversions of $\sigma$.

  By considering the natural action of the symmetric group $\sg_n$  on $B_n$, Stanley obtains an equivariant version of (\ref{rankseleq}) in \cite[Theorem 4.3]{St50}.  He proves that, in the representation of $\sg_n$ on the top homology of the order complex $\Delta\overline{(B_n)_S}$ of the rank-selected subposet $\overline{(B_n)_S}:=(B_n)_S\setminus \{\hat 0,\hat 1\}$,  the multiplicity   of the irreducible representation of $\sg_n$ indexed by the partition $\lambda$ is equal to the number of standard Young tableaux of shape $\lambda$ and descent set $S$.

  In \cite{BjWe}, Bj\"orner and Welker studied certain poset constructions arising in commutative algebra.  They defined the  {\it Rees product} $P \ast Q$ of ranked posets $(P, \leq_P)$ and $(Q,\leq_Q)$ with respective rank functions $r_P$ and $r_Q$.  The elements of $P \ast Q$ are those $(p,q) \in P \times Q$ satisfying $r_P(p) \geq r_Q(q)$.  A partial order on $P \ast Q$ is defined by setting $(p,q) \leq (p^\prime,q^\prime)$ if \begin{itemize} \item $p \leq_P p^\prime$, \item $q \leq_Q q^\prime$, and \item $r_P(p^\prime)-r_P(p) \geq r_Q(q^\prime)-r_Q(q)$. \end{itemize}
 
 Bj\"orner and Welker \cite{BjWe} proved that $P*Q$ is Cohen-Macaulay if $P$ is Cohen-Macaulay and $Q$ is acyclic and Cohen-Macaulay.   They  applied this result to the poset  $$R_n:=((B_n \setminus \{\emptyset\}) \ast C_n)\uplus\{\hat 0,\hat 1\},$$ where $C_n$ is the chain $\{0<1<\ldots <n-1\}$.      Bj\"orner and Welker  conjectured (and Jonsson proved in  \cite{Jo}) that 
 \begin{equation} \label{derangeeq} \mu_{R_n}(\hat 0, \hat 1) = (-1)^{n-1} d_n,\end{equation} 
 where 
$d_n$ is the number of derangements in the symmetric group $\mathfrak S_n$.  These results have been extended in several directions; see \cite{ShWa2}, \cite{LiShWa} and \cite{MuRe}.  Here we discuss just a few of these directions.

 Let $a_{n,j}$ be the number of permutations in $\sg_n$ with $j$ descents (or equivalently $j$ excedances; see (\ref{euler})).  Recall that the numbers $a_{n,j}$ are called {\em Eulerian numbers}.

\begin{thm}[Shareshian and Wachs {\cite[Theorem 1.2]{ShWa2}}] \label{rees1}
If $S \in B_n$ has size $m>0$, then $$\mu_{R_n}(\hat{0},(S,j))=(-1)^ma_{m,j}$$ for all $j \in \{0,1,\dots,m-1\}$.
\end{thm}

The conjecture of Bj\"orner and Welker (equation (\ref{derangeeq})) follows quickly from the standard inclusion-exclusion expression for $d_n$.  There are various ways to prove Theorem~\ref{rees1}.  In \cite{ShWa2} we apply the recursive definition of the M\"obius function 
to  a closely related Rees product poset
to obtain  the  exponential generating function formula, 
\begin{equation}\label{mobeq} 1+\sum_{n \geq 1} \sum_{j=0}^{n-1} (-1)^n \mu_{R_n}(\hat 0,( [n],j)) t^j \frac{z^n}{n!} = \frac{1-t}{e^{z(t-1)}-t}.\end{equation}
The result thus follows from Euler's  exponential generating function formula for the Eulerian numbers, given below in (\ref{euler}), and the fact that each lower interval $[\hat{0},(S,j)]$ of $R_n$ is isomorphic to the interval $[\hat 0, ([m],j)]$ in $R_{m}$.
In our joint paper \cite{LiShWa} with Linusson we gave another proof of Theorem~\ref{rees1}, which involves counting the ascent-free chains in an EL-labeling.  This enabled us to prove   generalizations of both Theorem~\ref{rees1} and  the Bj\"orner-Welker conjecture in which $B_n$ is replaced with an arbitrary product of chains and $\sym_n$ is replaced with a corresponding set of multiset permutations.

Inspired by  (\ref{qrankseleq}), we decided to consider next the $q$-analogue of $R_n$ obtained by replacing $B_n$ by $B_n(q)$.  
Let $$R_n(q):=(B_n(q) \setminus \{0\}) \ast C_n)\uplus\{\hat 0,\hat 1\}.$$  Since $B_n(q)$ is Cohen-Macaulay and $C_n$ is acyclic and Cohen-Macaulay, it follows from the general result of Bj\"orner and Welker that $R_n(q)$ is Cohen-Macaulay.  From our joint paper with Linusson \cite{LiShWa}, we get the stronger result  that $R_n(q)$ is EL-shellable.  
We have the following $q$-analog of Theorem~\ref{rees1}.

\begin{thm}[Shareshian and Wachs {\cite[Theorem~1.3]{ShWa2}}] \label{rees2} 
If $W \in B_n(q)$ has dimension $m>0$ and $0\le j \le m-1$, then 
$$\mu_{R_n(q)}(\hat{0},(W,j))= (-1)^{m} \sum_{\substack{\sigma \in \sg_m\\ \exc(\sigma) = j}} q^{{{m} \choose {2}}-\maj(\sigma)+j},$$
where  $\exc(\sigma)$ is the excedance number and $\maj(\sigma)$ is the major index of $\sigma$.
\end{thm}

The proof of Theorem~\ref{rees2} was substantially more difficult than any of the proofs of Theorem~\ref{rees1}. 
  We were originally led to conjecture the formula of Theorem~\ref{rees2} by inspection of  data.  To prove our conjectured formula, we  first derived a $q$-analog of (\ref{mobeq}).  This   served to reduce our conjectured formula to the $q$-analog of Euler's exponential generating function formula given in (\ref{majexceq}). To our surprise this $q$-analog of Euler's exponential generating function formula was new and   was not easy to prove.

 As a consequence of Theorem~\ref{rees2} we obtain the following $q$-analog of (\ref{derangeeq}).
\begin{cor} [{\cite[Corollary~1.4]{ShWa2}}] \label{rees2cor}
Let ${\mathcal D}_n$ be the set of derangements in $\sg_n$.  Then
 $$\mu_{R_n(q)}(\hat 0,\hat 1)= (-1)^{n-1}\sum_{\sigma \in {\mathcal D}_n}q^{{{n} \choose {2}}-\maj(\sigma)+\exc(\sigma)}.$$ 
\end{cor}

Combining the natural action of the symmetric group $\sg_n$ on $B_n$ with the trivial action of $\sg_n$ on $C_n$, we obtain an action on $R_n$. This action yields a representation of $\sg_n$ on  the unique nontrivial homology group $\widetilde{H}_{n-1}(\Delta (\overline{R_n}))$ of the order complex of  $\overline{R_n}:= R_n\setminus\{\hat 0,\hat 1\}$ and on  the unique nontrivial homology group $\widetilde{H}_{n-2}(\Delta I_j(B_n))$ of the order complex of  $I_j(B_n):=\{x \in R_n:\hat 0 < x<([n],j)\}$ for each $j$. 
We describe these representations (or $\sg_n$-modules) in \cite{ShWa2}.  The most efficient description involves the use of the Frobenius characteristic $\ch$, which assigns to each  $\sg_n$-module a homogeneous symmetric function of degree $n$.  Recall that the elementary symmetric function $\ee_n$ is the sum of all degree $n$ squarefree monomials in infinitely many variables $x_1,x_2,\ldots$.  We define $$E(z):=\sum_{n \geq 0}\ee_nz^n$$ and, for all $n \ge 1$,   $$[n]_t:= \frac {t^n -1}{t-1} = 1+t+\dots + t^{n-1}.$$ 

\begin{thm}[Shareshian and Wachs {\cite[Theorem 1.5 and Corollary 1.6]{ShWa2}}] \label{rees3}  The following equalities hold:

\begin{eqnarray} 
\label{charform} 1+\sum_{n \geq 1} \sum_{j=0}^{n-1}\ch(\widetilde{H}_{n-2}(\Delta I_j(B_n)))t^jz^n & = & \frac{(1-t)E(z)}{E(tz)-tE(z)} \\ \nonumber & = & \frac{\sum_{n \geq 0}\ee_n z^n}{1-\sum_{n \geq 2}t[n-1]_t\ee_nz^n}
\end{eqnarray}
and 
\begin{equation}\label{charform2} 1+\sum_{n \geq 1} \ch(\widetilde H_{n-2}(\Delta \overline{R_n}))z^n = \frac{1}{1-\sum_{n \geq 2}(n-1)\ee_n z^n}.\end{equation}
\end{thm}

Equation~(\ref{charform})  specializes to (\ref{mobeq}) when one applies a variant of  exponential specialization $\expon$  that takes $\ee_n$ to $\frac 1 {n!}$ 
for all $n\ge 0$ (cf. \cite[Sec. 7.8]{StEC2}).  
 Indeed, it is known  that if $V$ is an $\sg_n$-module then $\expon(\ch(V)) = \frac 1 {n!}\dim V$. 
Hence (\ref{charform})  specializes to 
$$1+\sum_{n \geq 1} \sum_{j=0}^{n-1}\dim \widetilde{H}_{n-2}(\Delta I_j(B_n))\,t^j \frac{z^n}{n!} = \frac{(1-t)e^z}{e^{tz}-te^z}.$$  
This is equivalent to (\ref{mobeq}) by the classical result of P. Hall relating M\"obius functions to reduced Euler characteristics (see \cite[Proposition 3.8.8]{StEC1}) and the Euler-Poincar\'e formula.  Similarly (\ref{charform2}) specializes to  the Bj\"orner-Welker conjecture (equation (\ref{derangeeq})).  

To prove Theorem~\ref{rees3} we use a technique of  Sundaram \cite{su}, which can be viewed as an equivariant version of the technique of using the recursive definition of M\"obius function to compute the Euler characteristic of an order complex.   An equivariant version of Theorem~\ref{rees1} is given in  Corollary~\ref{equireescor}.

 \section{$q$-Eulerian polynomials} \label{q-EulerSec}
 
\subsection{Euler's exponential generating function}  \label{eulersubsec} 
It follows from work of Euler in \cite{Eu}, an equidistribution result of
 MacMahon in \cite{MacM} and an observation of Riordan in \cite{Ri} that
\begin{equation} \label{euler}
1+ \sum_{n \geq 1}\sum_{\sigma \in \sg_n} t^{\des(\sigma) } \frac{z^n}{n!}=1+\sum_{n \geq 1}\sum_{\sigma \in \sg_n} t^{\exc(\sigma)} \frac{z^n}{n!}=\frac{1-t}{e^{z(t-1)}-t},
\end{equation}
where $\des(\sigma)$ and $\exc(\sigma)$ are the descent number and the excedance number of $\sigma$, respectively.  As referred to above, (\ref{euler}) is  {\it Euler's exponential generating function formula for the Eulerian polynomials}. 
 
  As mentioned in Section~\ref{posetsec}, in order to prove  Theorem~\ref{rees2} we needed to prove a certain $q$-analog of (\ref{euler}).  To obtain a $q$-analog of the Eulerian polynomials, one can combine a Mahonian statistic such as the major index $\maj$ or the inversion number $\inv$ with an Eulerian statistic such as $\exc$ or  $\des$. Such $q$-analogs have received considerable attention over the years.   A q-analog, involving $\maj$ and $\des$, of Euler's original definition of the Eulerian polynomials  was first obtained by MacMahon in  \cite[Vol. 2, Sect. IX, Ch. IV, Art. 462]{MacM} and rediscovered by Carlitz  in \cite{Ca}.   It follows from (\ref{qrankseleq}) that a $q$-analog, involving $\inv$ and $\des$, of the Eulerian polynomials  can be obtained from the rank-selected M\"obius invariant of  $B_n(q)$.  Using his theory of binomial posets,  Stanley proved in   \cite{St28}, a beautiful $q$-analog of (\ref{euler}) also involving involving $\inv$ and $\des$.  For $n \in \pp$, let $[n]_q!:=[n]_q [n-1]_q \cdots [1]_q$.
    
    \begin{thm}[Stanley  {\cite[p. 351]{St28}}] \label{invdes} Let $A_n^{\inv, \des}$ be the $q$-Eulerian polynomial defined by 
  $$A_n^{\inv,\des}(q,t) := \sum_{\sigma \in \sg_n} q^{\inv(\sigma)} t^{\des(\sigma)} .$$
   Then
\begin{equation} \label{staneq} 1+\sum_{n \geq 1}A_n^{\inv,\des}(q,t)\frac{z^n}{[n]!_q}=\frac{1-t}{\Exp_q(z(t-1))-t},  \end{equation}
where $$
\Exp_q(z):=\sum_{n \geq 0}q^{{n} \choose {2}}\frac{z^n}{[n]_q!}.
$$
\end{thm}

Our $q$-analog of (\ref{euler}) involves $\maj$ and $\exc$.
\begin{thm}[Shareshian and Wachs {\cite[Theorem 1.1]{ShWa1}, \cite[Corollary 1.3]{ShWa3}}] \label{majexc}
Let $A_n^{\maj, \exc}$ be the $q$-Eulerian polynomial defined by 
  $$A_n^{\maj,\exc}(q,t) := \sum_{\sigma \in \sg_n} q^{\maj(\sigma)} t^{\exc(\sigma)} .$$
   Then
\begin{equation} \label{majexceq}
1+\sum_{n \geq 1}A_n^{\maj,\exc}(q,t)\frac{z^n}{[n]!_q}=\frac{(1-tq)\exp_q(z)}{\exp_q(ztq)-tq\exp_q(z)},
 \end{equation}
 where
 $$
\exp_q(z):=\sum_{n \geq 0}\frac{z^n}{[n]_q!}.
$$
\end{thm}

From here on,   by  $q$-Eulerian polynomial  we mean 
$$A_n(q,t) := A_n^{\maj,\exc}(q,tq^{-1}) = \sum_{\sigma \in \sg_n} q^{\maj(\sigma)-\exc(\sigma)} t^{\exc(\sigma)},$$
and by $q$-Eulerian number  we mean the coefficients of $t^j$ in $A_n(q,t)$, that is
$$a_{n,j}(q) := \sum_{\substack{\sigma \in \sg_n\\ \exc(\sigma) = j}} q^{\maj(\sigma)-j} .$$
Now equation (\ref{majexceq}) can be rewritten as
\begin{equation} \label{majexceq2}
1+\sum_{n \geq 1}A_n(q,t)\frac{z^n}{[n]!_q}=\frac{(1-t)\exp_q(z)}{\exp_q(zt)-t\exp_q(z)}.
 \end{equation}

We prove (\ref{majexceq2}) by proving an identity of (quasi)symmetric functions and applying stable principal specialization.  (Recall that the stable principal specialization $\ps (f)$ of a quasisymmetric function $f(x_1,x_2,\ldots)$ is the power series in $q$ obtained by replacing each $x_i$ with $q^{i-1}$.)  As usual, $\h_n$ will denote the complete homogeneous symmetric function of degree $n$, that is, the sum of all degree $n$ monomials in $x_1,x_2,\ldots$.  Define $$H(z):=\sum_{n \geq 0}\h_nz^n.$$

For $S \subseteq [n-1]$, Gessel's fundamental quasisymmetric function $F_{n,S}$ is defined as 
\begin{equation} \label{gesdef} F_{n,S}:=\sum_{\substack{{i_1 \geq \ldots \geq i_n} \\ {i_j>i_{j+1} \forall j \in S}}}x_{i_1}x_{i_2} \cdots x_{i_n}.\end{equation}  
It is known (see \cite[Lemma 7.19.10]{StEC2}) that
\begin{equation} \label{psq}
\ps(F_{n,S})=\frac{q^{\sum_{i \in S}i}}{(q;q)_n},
\end{equation}
where $$(q;q)_n := \prod_{j=1}^{n}(1-q^j).$$
As $F_{n,\emptyset}=\h_n$, it follows that 
\begin{equation} \label{hzeq}
\ps(H(z(1-q)))=\exp_q(z).
\end{equation}

We define the alphabet $A:=[n] \cup [\overline{n}]$, where $[\overline{n}]:=\{\overline{i}:i \in [n]\}$.  From a permutation $\sigma=\sigma_1\ldots\sigma_n$, written in one-line notation, we obtain a word $\overline{\sigma}=\alpha_1\ldots\alpha_n$ in $A$ by replacing $\sigma_i$ with $\overline{\sigma_i}$ whenever $\sigma$ has an excedance at position $i$.  We order $A$ by $$\overline{1}<\ldots<\overline{n}<1<\ldots<n.$$ Now we define $\Dex(\sigma)$ to be the dsecent set  of $\overline{\sigma}$ in the given order on $A$.  A key point is that \begin{equation} \label{dex} \sum_{i \in \Dex(\sigma)}i=\maj(\sigma)-\exc(\sigma).\end{equation}
Now, for positive integers $n,j$, define the {\em Eulerian quasisymmetric function} $$Q_{n,j}:=\sum_{\substack {\sigma \in \sg_n \\ \exc(\sigma) = j} } F_{n,\Dex(\sigma)}.$$ 
It turns out that the quasisymmetric functions $Q_{n,j}$ are in fact symmetric functions. 
It follows from (\ref{psq}) and (\ref{dex}) that 
\begin{equation} \label{psQeq} \ps(Q_{n,j}) = \frac {a_{n,j}(q)}{(q;q)_n}.\end{equation}  
Upon replacing $z$ with $z(1-q)$, we obtain (\ref{majexceq2}) from (\ref{hzeq}), (\ref{psQeq}) and stable principal specialization of both sides of (\ref{qsid}) below.

\begin{thm}[Shareshian and Wachs {\cite[Theorem 1.2]{ShWa3}}]  \label{qsme} The following equalities hold:
\begin{eqnarray} \label{qsid}
1+\sum_{n \geq 1}\sum_{j=0}^{n-1}Q_{n,j}t^jz^n &=& \frac{(1-t)H(z)}{H(zt)-tH(z)} 
\\ \nonumber & = & \frac{H(z)}{1-\sum_{n \geq 2}t[n-1]_t\h_nz^n}.
\end{eqnarray}
\end{thm}

Let $\omega$ be the involution on the ring of symmetric functions that takes $\h_n$ to $\ee_n$.  By comparing Theorems~\ref{rees3} and~\ref{qsme} we obtain the following equivariant version of Theorem~\ref{rees1}.

\begin{cor} \label{equireescor} For all $j = 0,1,\dots,n-1$, 
$$\ch \widetilde{H}_{n-2}(\Delta I_j(B_n)) = \omega Q_{n,j}.$$
\end{cor}

 Equations~(\ref{majexceq2}) and~(\ref{qsid}) have been extended and applied in various ways; see e.g. \cite{ChGr, FoHa1,FoHa2,FoHa3,HaLiZe,HeWa,Hy, Li, LiShWa,SaShWa,ShWa3,ShWa4,ShWa5}.  In the remaining subsections  we discuss   just a few  of these applications and related results.

\subsection{Palindromicity and unimodality}  
Stanley proved many results about unimodality of sequences arising naturally in combinatorics.  Some of his proofs involved striking applications of algebraic geometry  (\cite{St42, St46}), commutative algebra (\cite{St88,St89}), convex geometry (\cite{St48}), and representations of groups, Lie algebras and Lie superalgebras (\cite{St41,St50,St57,St62,St84}). Some of this work is surveyed in \cite{St58,St72}).  
Motivated by Stanley's work on unimodality, we obtained in \cite{ShWa3} a $q$-analog of the well-known fact  that the 
Eulerian polynomials are palindromic and unimodal; see Theorem~\ref{unith} below.  

 Let $R$ be a $\Q$-algebra with basis ${\bold b}$.  For instance, $R$ can be the $\Q$-algebra of symmetric functions $\Lambda_\Q$ with  basis ${\bold h}:=\{\h_\lambda\}$ of complete homogeneous  symmetric functions or $R$ can be the $\Q$-algebra of polynomials $\Q[q]$ with basis $\q:=\{1,q,q^2,\dots\}$.  For 
$f,g \in R$, we say that $f \ge_{\bold b} g$ if the expansion of $f-g$ in the basis ${\bold b}$ has nonnegative 
coefficients.  We say that $f \in R$ is ${\bold b}$-positive if $f \ge_{\bold b} 0$ and that a polynomial $f(t) = \sum_{i=0}^{n} a_i t^i$ in $R[t]$ is ${\bold b}$-positive, ${\bold b}$-unimodal and 
palindromic with center of symmetry $\frac n 2$ if 
$$ 0 \le_{\bold b} a_{0} \le_{\bold b} a_{1} \le  \dots \le_{\bold b} a_{\lfloor {n\over 2} \rfloor} = a_{\lfloor {n+1\over 2} \rfloor} \ge_{\bold b} \dots \ge_{\bold b} a_{n-1}\ge_{\bold b} a_{n} \ge_{\bold b} 0,$$
and $a_i = a_{n-i}$ for all $i=0,\dots,n$.  

\begin{thm}[Shareshian and Wachs \cite{ShWa3}] \label{unith} The polynomial $\sum_{j=0}^{n-1} Q_{n,j}t^j \in \Lambda_\Q[t]$ is ${\bold h}$-positive, $\bold h$-unimodal and palindromic with center of symmetry $\frac {n-1}{2}$, and the polynomial $A_n(q,t):=\sum_{j=0}^{n-1} a_{n,j}(q)t^j \in \Q[q][t]$ is $\q$-positive, $\q$-unimodal and palindromic with center of symmetry $\frac {n-1}{2}$.
\end{thm}

There are various ways to prove this theorem.  In  his survey paper  \cite{St72}, Stanley presents various tools for establishing unimodality of sequences.   The most elementary tool  \cite[Proposition 1]{St72} states that the product of positive, unimodal, palindromic polynomials is  positive, unimodal, and palindromic.   Here we  describe how this  tool was used to prove Theorem~\ref{unith} in \cite[Appendix C]{ShWa4}.   In the next subsection we will see that the least elementary tool of \cite{St72} can also be used to prove the theorem.

Let $\left[ \begin{array}{c} n \\ k_1,  k_2, \dots, k_m \end{array} \right]_q$ denote the $q$-multinomial coefficient $\frac {[n]_q!} { [k_1]_q!\cdots [k_m]_q!}$.  In \cite[Appendix C]{ShWa4} we derive, by manipulating the right side of (\ref{majexceq2}), the following formula, from which Theorem~\ref{unith} is evident: \begin{equation}\label{newformq} A_n(q,t)= \sum_{m=1}^{\lfloor {n+1 \over 2} \rfloor}   \sum_{\substack{
k_1,\dots, k_m \ge 2 \\ \sum k_i = n+1}}\,\,\, \left[ \begin{array}{c} n \\ k_1-1,  k_2, \dots, k_m \end{array} \right]_q \,\,\,t^{m-1} \prod_{i=1}^m  [k_i-1]_{t} .\end{equation}
Indeed,   each  polynomial $t^{m-1} \prod_{i=1}^m  [k_i-1]_{t}$ is positive, unimodal and palindromic since each factor is.  The center of symmetry of each product is the sum of the centers of symmetries its factors. This sum is  $\frac{n-1}2$ in every case.  By expressing $A_n(q,t)$  as a sum  of $\q$-positive, $\q$-unimodal, palindromic polynomials all with the same center of symmetry we can conclude   that  $A_n(q,t)$ is $\q$-positive, $\q$-unimodal, and palindromic.  The result for $Q_{n,j}$ is proved similarly in \cite[Appendix C]{ShWa4}.

\subsection{Geometric interpretation: the hard Lefschetz theorem} \label{geomsec}
 One of Stanley's many major contributions to combinatorics is the idea of using the hard Lefschetz theorem to solve combinatorial problems; see  \cite{St42,St43,St46,St54,St58,St72}.  In this section we will discuss how the hard Lefschetz theorem can be used to prove Theorem~\ref{unith}.

It is well-known that the Eulerian polynomial $A_n(t):=A_n(1,t)$ equals the $h$-polynomial of the polar dual $\Delta_n$ of the $(n-1)$-dimensional permutohedron (or  equivalently, the type $A_{n-1}$ Coxeter complex). 

 Let $\Delta$ be a simplicial convex polytope of dimension $d$.  The $h$-polynomial $h_\Delta(t)$ of  $\Delta$ is defined by 
 $$h_\Delta(t) =\sum_{j=0}^d h_j t^{j} := \sum_{j=0}^{d} f_{d-1-j} (t-1)^{j},$$
where $f_i$ is the number of faces of $\Delta$ of dimension $i$.  The coefficient sequence $(h_0,h_1,\dots,h_d)$ is known as the $h$-vector of $\Delta$.  The Dehn-Sommerville equations state that  $h_\Delta(t)$ is palindromic.  Stanley's use of  the hard Lefschetz  theorem to prove the following theorem is a milestone in combinatorics.
\begin{thm}[Stanley \cite{St46}] \label{unihvecth}  If $\Delta$ be a simplicial convex polytope then $h_\Delta(t)$ is unimodal (as well as palindromic).
\end{thm}
Theorem~\ref{unihvecth} is an important consequence of the necessity part of the celebrated $g$-theorem.  (The sufficiency part was proved by Billera and Lee in \cite{BiLe}.)  We describe the role of the  hard Lefschetz theorem in the proof of Theorem~\ref{unihvecth}.  We may assume that the polytope $\Delta$ has rational vertices. There is a complex projective variety ${\mathcal V}_\Delta$, of (complex) dimension $d$, naturally associated to the rational polytope $\Delta$, namely, the  toric variety associated with $\Delta$.  In \cite{Da} Danilov proved that $$h_\Delta(t) = \sum_{i=0}^d \dim H^{2i}({\mathcal V}_\Delta) t^i,$$ where $H^{2i}$ denotes cohomology in degree $2i$.  Stanley saw that the hard Lefschetz theorem can be applied in this situation.  This theorem provides an injective map from $H^{2i-2}({\mathcal V}_\Delta)$ to $H^{2i}({\mathcal V}_\Delta)$ when $1 \le i \le \frac d 2$ (as well as an isomorphism from $H^{2i}({\mathcal V}_\Delta)$ to $H^{2d-2i}({\mathcal V}_\Delta)$). Thus the polynomial $\sum_{i=0}^d \dim H^{2i}({\mathcal V}_\Delta) t^i$ is palindromic and unimodal.
By setting $\Delta$ equal to the dual permutohedron $\Delta_n$, we obtain a geometric proof that the Eulerian polynomial $A_n(t)= h_{\Delta_n}(t)$ is palindromic and unimodal.  

 To obtain a geometric proof of {$\q$-unimodality and palidromicity of the $q$-Eulerian polynomials $A_n(q,t)$,  we use  our Eulerian quasisymmetric functions $Q_{n,j}$ and a result of Procesi and Stanley.
The reflection action of $\sg_n$ on $\Delta_n$ determines a linear representation of $\sg_n$ on each cohomology  of ${\mathcal V}_{\Delta_n}$.  
Stanley \cite[Proposition 12]{St72}  used a recurrence of Procesi \cite{Pr} to obtain the following generating function formula for the Frobenius characteristic of the $\sg_n$-module $H^{2j}  ({\mathcal V}_{\Delta_n})$:
\begin{equation} \label{prosta} 1+\sum_{n \geq 1} \sum_{j=1}^{n-1} \ch(H^{2j}({\mathcal V}_{\Delta_n}))t^jz^n=\frac{(1-t)H(z)}{H(zt)-tH(z)} . \end{equation}
By comparing (\ref{prosta}) with Theorem~\ref{qsme}, we conclude that for $0\le j \le n-1$,
\begin{equation} \label{Qtoreq} Q_{n,j} = \ch(H^{2j}({\mathcal V}_{\Delta_n})).\end{equation}
Now from (\ref{psQeq}) we obtain
$$ (q;q)_n\, \ps(\ch(H^{2j}({\mathcal V}_{\Delta_n})) )=a_{n,j}(q).$$
Next we use the observation of Stanley in \cite[p. 528]{St72} that since the hard Lefschetz map commutes with the action of $\sg_n$, 
the polynomial $\sum_{i=0}^{n-1} \ch(H^{2j}({\mathcal V}_{\Delta_n})) t^j \in \Lambda_\Q[t]$ is palindromic and Schur-unimodal.  Since $(q;q)_n\, \ps(s_\lambda) \in \N[q]$ for every Schur function $s_\lambda$ of degree $n$,  we conclude that $A_n(q,t)$ is palindromic and $\q$-unimodal.

We remark that results analogous to Stanley's formula (\ref{prosta}) for Coxeter complexes associated to arbitrary crystallographic root systems were obtained by Stembridge in \cite{Ste2} and by Dolgachev-Lunts in \cite{DoLu}.

\subsection{Refinements} The symmetric functions $Q_{n,j}$ and the $q$-Eulerian numbers $a_{n,j}(q)$ can be refined.  For a partition $\lambda$ of $n$ and a nonnegative integer $j$, we define 
${\sg}_{\lambda,j}$ to be the set of all $\sigma \in \sg_n$ with $j$ excedances and cycle type $\lambda$.  Now set
$$
Q_{\lambda,j}:=\sum_{\sigma \in {\sg}_{\lambda,j}} F_{n,\Dex(\sigma)}
$$
and $$a_{\lambda,j}(q):=\sum_{\sigma \in {\sg}_{\lambda,j}} q^{\maj(\sigma)-j}.$$
Note that by (\ref{psq}) and (\ref{dex}),
$$a_{\lambda,j}(q) =  (q;q)_n \, \ps(Q_{\lambda,j}).$$
It turns out that the quasisymmetric functions $Q_{\lambda,j}$ are symmetric and have some remarkable properties, see \cite{ShWa3, HeWa,SaShWa}. 

In \cite{Br}, Brenti proved that $A_\lambda(t):= \sum_{j\ge 0} a_{\lambda,j}(1) t^j$ is unimodal and palindromic. We conjectured the following symmetric function analog and $q$-analog of this result  in \cite{ShWa3}.
 \begin{thm}[Henderson and Wachs \cite{HeWa}] The polynomial $\sum_{j\ge 0} Q_{\lambda,j} t^j$ is Schur-positive, Schur-unimodal, and palindromic.  Consequently, the polynomial $A_\lambda(q,t):=\sum_{j\ge 0} a_{\lambda,j}(q) t^j$ is $\q$-unimodal, and palindromic.  
 \end{thm} 
We do not have a geometric interpretation of this result like that for $\sum_{j=0}^{n-1} Q_{n,j} t^j$ and $A_n(q,t)$ given in Section~\ref{geomsec}.  Schur-positivity is proved by deriving a plethystic formula for $Q_{\lambda,j}$ from one obtained by us in \cite{ShWa3}.   Schur-unimodality is established by using the plethystic formula to construct  an $\sg_n$-module $V_{\lambda,j}$ whose Frobenius characteristic is $Q_{\lambda,j}$ and an $\sg_n$-module monomorphism from $V_{\lambda,j-1}$ to  $V_{\lambda,j}$ when $1 \le j \le \frac {n-k} 2$, where $k$ is the number of parts of size $1$ in $\lambda$.  Many other related unimodality results  and conjectures appear in \cite{ShWa3,HeWa}.
 
 The $Q_{\lambda,j}$ also play a critical role  in our joint paper \cite{SaShWa} with Sagan, in which we show that the polynomials $a_{\lambda,j}(q)$ appear in an instance of the cyclic sieving phenomenon of Reiner, Stanton, and White (see \cite{ReStWh}).

\begin{thm}[Sagan, Shareshian, and Wachs {\cite[Theorem 1.2]{SaShWa}}] \label{csp}
Let $\gamma_n \in \sg_n$ be the $n$-cycle $(1,2,\ldots,n)$.  Then, for each partition $\lambda$ of $n$ and each $j \in \{0,\ldots,n-1\}$, the subgroup $\langle \gamma_n \rangle$ generated by $\gamma_n$ acts by conjugation on ${\mathcal S}_{\lambda,j}$.   If $\tau \in \langle \gamma_n \rangle$ has order $d$, then the fixed point set of $\tau$ in this action has size $a_{\lambda,j}(e^{2\pi i/d})$.
\end{thm}

\section{Chromatic quasisymmetric functions} \label{csfsec}

\subsection{An email message from Stanley}  After distributing a preliminary version of our research announcement  \cite{ShWa1} on $q$-Eulerian polynomials, we received an email message from Stanley. Therein, he pointed out to us that   a characterization of the Eulerian quasisymmetric functions  $Q_{n,j}$ given  in \cite[Proposition 2.4]{ShWa1} (see also \cite[Theorem 3.6]{ShWa3}) is equivalent to the characterization
\begin{equation}\label{smirnoveq} \omega Q_{n,j} = \sum_{\substack {w \in W_{n}\\ \des(w) =j}} x_w. \end{equation}
Here  \begin{itemize} 
\item  $W_{n,}$ is the set of all words $w:=w_1\ldots w_n$ over the alphabet of positive integers $\pp$ satisfying $w_i \neq w_{i+1}$ for all $i \in [n-1]$,
\item $\des(w)$ is the number of descents of $w$, and
\item $x_w:= x_{w_1}x_{w_2} \cdots x_{w_n}$.  
\end{itemize}
Indeed, this equivalence follows from $P$-partition reciprocity, which Stanley introduced in \cite{St9}; see  \cite[equation (7.7)]{ShWa3}.

Theorem~\ref{qsme} thus gives a generating function formula for $ \sum_{w \in W_{n}} x_w t^{\des(w)}$, which refines the following formula of 
 Carlitz, Scoville and Vaughan (see \cite{CaScVa}), 
\begin{equation} \label{csveq}
1+\sum_{n \geq 1}\sum_{w \in W_n} x_w z^n=\frac{\sum_{n \geq 0} \ee_nz^n}{1-\sum_{n \geq 2}(n-1)\ee_nz^n}.
\end{equation}
Stanley observes in \cite[Proposition 5.3]{St100} that since the words in $W_n$ can be viewed as proper colorings of the path 
graph $G_n:= ([n], \{ \{i,i+1\} : i \in [n-1]\}$, the symmetric function $\sum_{w \in W_n} x_w $ can be viewed as the 
chromatic symmetric function of $G_n$. 
These observations of Stanley brought us from the world of Eulerian polynomials to the world of chromatic symmetric functions.  We consider in \cite{ShWa4,ShWa5} a refined notion of chromatic symmetric 
function, one for which $\sum_{w \in W_n} x_w t^{\des(w)}$ is the refined chromatic symmetric function for
$G_n$.   Before describing our refinement we review chromatic symmetric functions.

\subsection{Stanley's chromatic symmetric function}
Let $G=(V,E)$ be a simple finite graph with no loops.  A proper $\pp$-coloring of $G$ is a function $\kappa:V\to \pp$ such that $\kappa(u) \neq \kappa(v)$ whenever $\{u,v\}$ lies in the edge set $E$.  For such a coloring $\kappa$, set 
$$
x_\kappa :=\prod_{v \in V}x_{\kappa(v)}.
$$
 In \cite{St100}, Stanley defined the chromatic symmetric function $$X_G(\x):=\sum_\kappa x_\kappa,$$ where the sum is taken over all proper $\pp$-colorings $\kappa$ of $G$.  It is apparent that $X_G$ is indeed a symmetric function and that by setting $x_i = 1$ for all $i \in [m]$ and $x_i = 0$ for all $i \ge m+1$, one gets   the chromatic polynomial $\chi_G$ of $G$ evaluated at $m$.  Stanley studied various aspects of $X_G$ in \cite{St100} and \cite{St101}.  In particular, he examined expansions of $X_G$ in various well-studied bases for the ring of symmetric functions.  
 
 Among the many interesting questions and theorems about $X_G$ are ones that arise when the structure of $G$ is restricted appropriately.  Given a finite poset $P$, the incomparability graph $\inc(P)$ has vertex set $P$.  Its edge set consists of all $\{p,q\}$ such that neither $p<q$ nor $q<p$ holds in $P$.  For positive integers $a,b$, we call $P$ $(a+b)$-free if there is no induced subposet of $P$ that is the disjoint union of a chain of $a$ elements and a chain of $b$ elements.
  The next conjecture was first stated in the form given below in \cite{St100}.  It is a generalization of a particular case of a conjecture of Stembridge on immanants, see \cite{Ste1}.  The transformation of this case to a statement about chromatic symmetric functions is achieved implicitly by Stanley and Stembridge in \cite{StSt}.

\begin{con}[Stanley/Stembridge, see {\cite[Conjecture 5.1]{St100}}] \label{ststcon}
If $P$ is a $(3+1)$-free poset then $X_{\inc(P)}$ is $\ee$-positive.
\end{con}

Stanley and Stembridge showed in \cite{StSt} that Conjecture~\ref{ststcon} is true for all posets in some interesting infinite classes.  For example, if $P$ has no chain of length three then $X_{\inc(P)}$ is $\ee$-positive.  Gasharov showed in \cite{Ga1} that $X_{\inc(P)}$ is Schur-positive when $P$ is $(3+ 1)$-free by describing the coefficients of the Schur functions in the  expansion of $X_{\inc(P)}$ in terms of tableaux.  Guay-Paquet showed in \cite{Gu-Pa} that if Conjecture~\ref{ststcon} holds for all $P$ that are both $(3+1)$-free and $(2+2)$-free, then it holds in general.  Other work on chromatic symmetric functions can be found in \cite{Ga2, Ch1, Ch, Wo, NoWe, MaMoWa, Hu}.

\subsection{A quasisymmetric refinement} \label{chromquasisec}
Given a proper $\pp$-coloring $\kappa$ of a graph $G=([n],E)$, we define the ascent number $\asc(\kappa)$ to be the number of edges $\{i,j\} \in E$ such that $i<j$ and $\kappa(i)<\kappa(j)$. In \cite{ShWa4,ShWa5}, we define the {\it chromatic quasisymmetric function} 
$$X_G(\x,t):=\sum_\kappa t^{\asc(\kappa)}x_\kappa,$$ where the sum is taken over all proper $\pp$-colorings $\kappa$.  So, $X_G(\x,t)$ refines the chromatic symmetric function $X_G(\x)$.

\begin{example}  \label{pathex} We show in \cite[Example 3.2]{ShWa5}  that if $G$ is the path $1-2-3$ then $$X_{G}(\x,t) = \ee_3 + (\ee_3 + \ee_{2,1})t  + \ee_3t^2$$ and if $G$ is the path $1-3-2$ then $$X_{G}(\x,t) = (\ee_3+F_{3,\{1\}}) + 2 \ee_3 t + (\ee_3 + F_{3,\{2\}})t^2,$$ where $F_{n,S}$ denotes Gessel's fundamental quasisymmetric function defined in (\ref{gesdef}).
\end{example}
Upon considering the paths with three vertices given in Example~\ref{pathex}, one observes two clear differences between $X_G(\x)$ and $X_G(\x,t)$.  First, $X_G(\x,t)$ depends not only on the isomorphism type of $G$ but also on the labeling of the vertices of $G$ with the elements of $[n]$.  Second, if we consider $X_G(\x,t)$ to lie in the polynomial ring in the variable $t$ with coefficients in the ring of power series in $x_1,x_2,\ldots$, then the coefficient of $t^j$ in $X_G(\x,t)$ need not be a symmetric function.  This coefficient is, however, a quasisymmetric function, hence the name.

In \cite{ShWa4,ShWa5} we consider a class of graphs on $[n]$  for which the coefficient of each $t^j$ in $X_G(\x,t)$ is a symmetric function.    
 Choose a finite set of closed intervals $[a_i,a_i+1]$ ($1 \leq i \leq n$) of length one on the real line, with $a_i<a_{i+1}$ for $1 \leq i \leq n-1$.  The associated natural unit interval order $P$ is the poset on $[n]$ in which $i<_P j$ if $a_i+1<a_j$.  A natural unit interval order is both $(3+1)$-free and $(2+2)$-free.  Conversely, every finite poset that is both $(3+1)$-free and $(2+2)$-free is isomorphic with a unique natural unit interval order.  We show  that if $G$ is the incomparability graph of a natural unit interval order then  $X_G(\x,t)$ is  symmetric  in $\x$ and   palindromic as a polynomial in $t$.   

\begin{con}[Shareshian and Wachs {\cite[Conjecture 5.1]{ShWa5}}]  \label{quasistan} Let $G$ be the incomparability graph of a natural unit interval order.   Then  the palindromic polynomial $X_{G}(\x,t)$  is $\ee$-positive and $\ee$-unimodal.   \end{con} 

Conjecture~\ref{quasistan} refines Conjecture~\ref{ststcon} for posets that are not only $(3+1)$-free but also $(2+2)$-free.  However, given the result of Guay-Paquet mentioned above, Conjecture~\ref{quasistan} implies Conjecture~\ref{ststcon}.

\begin{example}  The path $G_n:=1-2-\cdots - n$ is the incomparability graph of a natural unit interval order.  In \cite{ShWa5} we use  (\ref{smirnoveq}) and (\ref{qsid}) to derive
$$X_{G_n}(\x,t) =  \sum_{m=1}^{\lfloor {n+1 \over 2} \rfloor}   \sum_{\substack{
k_1,\dots, k_m \ge 2 \\ \sum k_i = n+1}}\,\,\, \ee_{ (k_1-1,  k_2, \dots, k_m)} \,\,\,t^{m-1} \prod_{i=1}^m  [k_i-1]_{t} \,\,.$$  As discussed in \cite[Appendix C]{ShWa5}, by expressing $X_{G_n}(\x,t)$ in the formula above  as a sum of $\ee$-positive, $\ee$-unimodal, palindromic polynomials with the same center of symmetry,   one can conclude that the conjecture holds for $G_n$.
\end{example}

Other classes of examples for which Conjecture~\ref{quasistan} holds are given in \cite[Section 8]{ShWa5} and weaker versions of the conjecture are also discussed in \cite{ShWa5}.  For instance, in \cite[Section 6]{ShWa5} we establish Schur-positivity of $X_G(\x,t)$ when $G$ is the incomparability graph of a natural unit interval order.  In fact, we give a formula for the coefficient of each Schur function $s_\lambda$  in the expansion of  $X_G(\x,t)$, refining the unit interval order case of the result of Gasharov in \cite{Ga1} mentioned above.  Given a poset $P$ on $[n]$ and a partition $\lambda$ of $n$, Gasharov defines a $P$-tableau
  of shape $\lambda$} to be  a filling  of a Young diagram of shape $\lambda$ (in English notation)
with elements of $P$ such that
\begin{itemize}
  \item each element of $P$ appears exactly once,
  \item if $y \in P$ appears immediately to the right of $x \in P$ then $y>_P x$, 
  \item if $y \in P$  appears immediately below  $x\in P$ then $y \not <_P x$.
\end{itemize}
Let ${\mathcal T}_P$ be the set of all $P$-tableaux. 
 For  $T \in \mathcal T_P$ and graph $G=([n],E)$, define a {\em $G$-inversion} of $T$ to be  an edge $\{i,j\} \in E$ such that 
 $i<j$ and $i$ appears below  $j$ in $T$ (not necessarily in the same column).  Let $\inv_G(T)$ be the number of $G$-inversions of $T$
  and let  $\lambda(T)$ be the shape of $T$.
   \begin{example} \label{tabex} Let  $P:=P_{n,r}$ be the poset on $[n]$ with $i <_P j$ if $j-i\ge r$ and let $G_{n,r}$ be the  incomparability graph  of $P_{n,r}$.  That is,  $G_{n,r} = ([n],\{\{i,j\} : i,j \in [n],  0<|j-i| < r\}$.  If $P:=P_{9,3}$ and  $G:=G_{9,3}$ then $$T=\tableau{ {2} & {6}  & {9}    \\ {1} & {4} & {8} \\  {3} & {7}  \\ {5} }$$
 is a $P$-tableau of shape $(3,3,2,1)$ and 
$$\inv_G(T) = |\{ \{1,2\},  \{3,4\}, \{4,6\}, \{5,6\},\{5,7\},  \{7,8\} , \{7,9\},\{8,9\} \}| = 8 .$$
\end{example}

Upon setting $t=1$ in (\ref{gengash}) below, we obtain a formula that was shown to hold for all $(3+1)$-free posets by Gasharov in \cite{Ga1}.
\begin{thm}[Shareshian and Wachs {\cite[Theorem 6.3]{ShWa5}}] \label{schurcon}
Let $G$ be the incomparability graph of  be a natural unit interval order $P$.  Then
\begin{equation} \label{gengash}
X_G(\x,t)=\sum_{T \in {\mathcal T}_P}t^{\inv_G(T)}s_{\lambda(T)}.
\end{equation}
Consequently, $X_G(\x,t)$ is Schur-positive.
\end{thm}

We use Theorem~\ref{schurcon} and the Murnaghan-Nakayama rule to prove in \cite[Section 7]{ShWa5} that for the incomparability graph $G:=([n],E)$ of a natural unit interval order,  the coefficient of  $\p_n$ in the power-sum symmetric function basis expansion of $\omega X_G(\x,t)$ is 
$$\frac {[n]_t} {n} \prod_{j=2}^n [b_j]_t,$$ where $b_j = |\{i,j\} \in E :  i < j\}$. Palindromicity and unimodality of the coefficient is a consequence of this formula.  We also make two (equivalent) conjectures describing the  coefficient of $\p_\lambda$ for  each  partition $\lambda$ of $n$, which can be shown to refine the corresponding result of Stanley for $X_{G}(\x)$ in \cite{St100}.  These conjectures have been proved by Athanasiadis in \cite{Ath}. 

Chromatic quasisymmetric functions of incomparability graphs of natural unit interval orders are also related to representations of type A Hecke algebras, as explained in \cite{ClHyShSk}.

\subsection{Generalized $q$-Eulerian polynomials and an exercise in Stanley's EC1}
It follows from (\ref{psQeq}) and (\ref{smirnoveq}) that the stable principal specialization of the chromatic quasisymmetric function of the path $G_n$ is given by
\begin{equation} \label{specpath} \ps(\omega X_{G_n}(\x,t)) = \ps(\sum_{j=0}^{n-1} Q_{n,j} t^j) = (q;q)_n^{-1} {A_n(q,t)}.\end{equation} 
To compute $\ps(\omega X_{G}(\x,t))$  for an arbitrary incomparability graph $G=\inc(P)$, we use a refinement of Chow's expansion of $X_{G}(\x)$ in the fundamental quasisymmetric function basis (see \cite{Ch}). In \cite[Section 3]{ShWa5}, we show that
\begin{equation} \label{chowquasi}
X_{G}(\x,t) = \sum_{\sigma \in \sg_n} F_{n,\Des_P(\sigma)} t^{\inv_G(\sigma)},\end{equation} 
where 
$$\Des_P(\sigma):= \{i \in [n-1]: \sigma(i) >_P \sigma(i+1) \}$$
and
$$\inv_G(\sigma) := |\{ \{\sigma(i),\sigma(j)\} \in E(G): i<j \in [n], \sigma(i) > \sigma(j)\}|.$$
For a certain class of natural unit interval orders, the stable principal specialization of the corresponding chromatic quasisymmetric function has a particularly attractive description, which we present here.

 For $\sigma \in \sg_n$, let
\begin{eqnarray*} \inv_{<r}(\sigma) &:=& |\{(i,j) : i < j, \,\,0<\sigma(i)-\sigma(j) < r\}|\\ 
\maj_{\ge r}(\sigma) &:=& \sum_{i : \sigma(i) - \sigma(i+1) \ge r}  i
.\end{eqnarray*}
The permutation statistic $\inv_{<r}+ \maj_{\ge r}$ was introduced by Rawlings  \cite{Ra} who proved that 
it is Mahonian for all $r$, that is
$$\sum_{\sigma \in \sg_n} q^{\maj_{\ge r}(\sigma)+ \inv_{<r}(\sigma)} = [n]_q !.$$
Note that this Mahonian statistic interpolates between $\maj$ (when $r=1$) and $\inv$  (when $r=n$).

Now consider the polynomial obtained by splitting the Rawlings statistic as follows 
$$ A^{(r)}_n(q,t) := \sum_{\sigma \in \sg_n} q^{\maj_{\ge r}(\sigma)} t^{\inv_{<r}(\sigma)}.$$
 As in Example~\ref{tabex}, for each $r \in [n]$,  let $P_{n,r}$  be the poset on  $[n]$ with order relation given by $i <_{P_{n,r}} j$ if $j-i \ge r$. 
Note that when $r=2$, 
$\inc(P_{n,r})$ is the path $G_n$. 
By  (\ref{chowquasi}) and (\ref{psq}),  we have
$$A^{(r)}_n(q,t) = (q;q)_n \Ps(\omega X_{\inc(P_{n,r})}(x,t)).
$$
It therefore follows from (\ref{specpath}) that 
\begin{equation} \label{genEulereq} A_n^{(2)}(q,t) = A_n(q,t).\end{equation} 
This justifies  calling $A_n^{(r)}(q,t)$  a generalized $q$-Eulerian polynomial. 
Since $\inv_{<2}(\sigma) = \des(\sigma^{-1})$, the $q=1$ case of (\ref{genEulereq}) is equivalent to the fact that $\des$ and $\exc$ are equidistributed on $\sg_n$, for which there is a well-known bijective proof.  In \cite[Problem 9.8]{ShWa5} we pose the problem of finding a direct bijective proof of (\ref{genEulereq}).

It is a consequence of the palindromicity of $X_{\inc(P_{n,r})}(\x,t)$ that the generalized $q$-Eulerian polynomials are palindromic as polynomials in $t$.  Conjecture~\ref{quasistan} implies the following conjecture.

\begin{con}[{\cite[Conjecture 9.6]{ShWa5}}] \label{genqconj} For all $r \in [n]$, the palindromic polynomial $A_n^{(r)}(q,t)$ is $\q$-unimodal.
\end{con}

After formulating this conjecture we discovered  that the $q=1$ case of the conjecture appears  as one of the 203 exercises of Chapter 1 of Stanley's EC1 (see \cite[Exercise 1.50 f]{StEC1}).  
The solution in \cite[page 157]{StEC1} is given by the following result, which involves a certain variety called a Hessenberg variety whose definition is given in the next subsection.

\begin{thm}[De Mari and Shayman \cite{DeSh}] Let $\mathcal H_{n,r}$ be the type $A_{n-1}$ regular semisimple Hessenberg variety of degree $r$.  Then 
$$A_n^{(r)}(1,t) = \sum_{j \ge 0} \dim H^{2j}(\mathcal H_{n,r}) t^j $$
Consequently by the hard Lefschetz theorem, $A_n^{(r)}(1,t)$ is palindromic and unimodal.
\end{thm}

Stanley asked for a more elementary proof of unimodality of $A^{(r)}_n(1,t)$ in \cite[page 157]{StEC1} .  (As far as we know this is still open.)  

The result of De Mari and Shayman suggested to us an approach to proving Conjecture~\ref{genqconj} along the lines of the geometric proof of $\q$-unimodality of $A_n(q,t)$ given in Section~\ref{geomsec}.  What was needed was a representation of the symmetric group on the cohomology of the Hessenberg variety $\mathcal H_{n,r}$ whose Frobenius characteristic is  $\omega X_{\inc(P_{n,r})}(\x,t)$.   
In the next subsection we discuss a promising candidate for such a representation.

\subsection{Hessenberg varieties}
A weakly increasing sequence ${\bf m}:=(m_1,\ldots,m_{n-1})$ of integers in $[n]$  is called a {\em Hessenberg vector} if  $m_i \ge i$ for each $i$.
 Let $\mathcal F_n$ be the variety of all  flags of subspaces  $$ F: V_1 \subset V_2 \subset \cdots \subset V_n = \C^n$$ with $\dim V_i = i$.  Fix $s \in GL_n(\C)$ such that $s$ is diagonal with $n$  distinct eigenvalues.   The {\em regular semisimple Hessenberg variety of type $A$} associated to the Hessenberg vector ${\bf m}$ is
\[
\hess ({\bf m}):=\{F \in \mathcal F_n \mid s V_i \subseteq V_{m_i}  \mbox{ for all } i \in [n-1]\}.
\]
The degree $r$ Hessenberg variety $\mathcal H_{n,r}:=\mathcal H(r,r+1,\dots,n,\dots,n)$ was studied initially by De Mari and Shayman in \cite{DeSh}.  Further work on regular semisimple Hessenberg varieties associated to arbitrary root systems was done by De Mari, Procesi and Shayman in \cite{DePrSh}.

Given a Hessenberg vector ${\bf m}$,  define $P({\bf m})$ to be the poset on  $[n]$ defined by $i <_{P({\bf m})} j$ if $ j > m_i$.  
 In \cite{ShWa5} we observe that  $P({\bf m})$ is a natural unit interval order  for all Hessenberg vectors ${\bf m}$, and  for every natural unit interval order  $P$ there is a unique Hessenberg vector ${\bf m}$ satisfying $P=P({\bf m})$.

As observed in \cite{DePrSh}, ${\mathcal H}_{n,2}:={\mathcal H}(2,3,\dots,n)$
 is the toric variety ${\mathcal V}_{\Delta_n}$.  This variety admits an $\sg_n$ action as described in Section~\ref{geomsec}.  We do not know whether for general Hessenberg vectors ${\bf m}$, each ${\mathcal H}({\bf m})$ admits a faithful $\sg_n$ action.  
However, there is a nice representation of $\sg_n$ on $H^\ast({\mathcal H}({\bf m}))$, as described by Tymoczko in \cite{Ty2}.  Let $T$ be the group of diagonal matrices in $GL_n(\cc)$.  It is straightforward to confirm that if $u \in T$ and the flag ${\mathcal F}$ lies in ${\mathcal H}({\bf m})$, then $u{\mathcal F} \in {\mathcal H}({\bf m})$.  This action of $T$ on ${\mathcal H}({\bf m})$ satisfies the technical conditions necessary to apply the theory of Goresky-Kottwitz-MacPherson on torus actions (see \cite{GoKoMacP}, and see \cite{Ty1} for a survey of GKM theory).  The GKM theory says that the cohomology of ${\mathcal H}({\bf m})$ is determined by the moment graph $M$ associated to the $T$-action.  It turns out that $M$ admits an $\sg_n$ action, which determines a representation of $\sg_n$ on each cohomology  of ${\mathcal H}({\bf m})$.  

When ${\bf m} = (2,3,\dots,n)$,  the representation of $\sg_n$ on $H^{2j}({\mathcal H}({\bf m}))$, for each $j$, is isomorphic to the representation of $\sg_n$ on $H^{2j} ({\mathcal V}_{\Delta_n})$ discussed in Section~\ref{geomsec}; see \cite{Ty2}.  It therefore follows from (\ref{Qtoreq}) and (\ref{smirnoveq}) that the following conjecture holds when ${\bf m} = (2,3,\dots,n)$.
\begin{con}[Shareshian and Wachs {\cite[Conjecture 10.1]{ShWa5}}] \label{hesscon} For all Hessenberg vectors ${\bf m}$, $$ \omega X_{\inc(P({\bf m}))}(\x,t) = \sum_{j\ge 0} \ch (H^{2j}({\mathcal H}({\bf m}))) t^j.$$
\end{con}

There is considerable evidence in favor of Conjecture~\ref{hesscon}, see  \cite[Section 5]{ShWa4} and \cite[Section 10]{ShWa5}. 
Since the hard Lefschetz map on $H^{*}({\mathcal H}({\bf m}))$ commutes with the action of the symmetric group,  the conjecture implies that $X_{\inc(P({\bf m}))}(\x;t)$ is Schur-unimodal, which in turn implies that Conjecture~\ref{genqconj} holds.  

We aim to carry out the following two-step process, thereby repaying Richard Stanley in some small measure for all of the wonderful ideas he has shared with us and our colleagues.
\begin{enumerate}
\item Prove Conjecture~\ref{hesscon}.
\item Understand the action of $\sg_n$ on $H^\ast({\mathcal H}({\bf m}))$ well enough to prove Conjecture~\ref{quasistan}, thereby proving the Stanley-Stembridge $\ee$-positivity conjecture (Conjecture~\ref{ststcon}).
\end{enumerate}

\bibliographystyle{amsplain}

\end{document}